\documentclass[12pt,leqno,draft]{article}

\newtheorem{theorem}{Theorem}
\newtheorem{lemma}[theorem]{Lemma}
\newtheorem{proposition}[theorem]{Proposition}
\newtheorem{sublemma}[theorem]{Sublemma}
\newtheorem{definition}[theorem]{Definition}
\newtheorem{corollary}[theorem]{Corollary}
\newtheorem{problem}[theorem]{Problem}
\newtheorem{remark}[theorem]{Remark}
\newtheorem{claim}[theorem]{Claim}
\newtheorem{assumptions}[theorem]{Assumptions}
\newtheorem{examples}[theorem]{Examples}
\newtheorem{question}[theorem]{Question}
\newtheorem{sassumptions}[theorem]{Standing Assumptions}
\newtheorem{sassumption}[theorem]{Standing Assumption}
\newtheorem{conjecture}[theorem]{Conjecture}

\newcommand{\begintheorem}{\addtocounter{equation}{1}\begin{theorem}}
\newcommand{\beginlemma}{\addtocounter{equation}{1}\begin{lemma}}
\newcommand{\beginproposition}{\addtocounter{equation}{1}\begin{proposition}}
\newcommand{\beginsublemma}{\addtocounter{equation}{1}\begin{sublemma}}
\newcommand{\begindefinition}{\addtocounter{equation}{1}\begin{definition}}
\newcommand{\begincorollary}{\addtocounter{equation}{1}\begin{corollary}}
\newcommand{\beginproblem}{\addtocounter{equation}{1}\begin{problem}}
\newcommand{\beginremark}{\addtocounter{equation}{1}\begin{remark}}
\newcommand{\beginclaim}{\addtocounter{equation}{1}\begin{claim}}
\newcommand{\beginassumptions}{\addtocounter{equation}{1}\begin{assumptions}}
\newcommand{\beginexamples}{\addtocounter{equation}{1}\begin{examples}}
\newcommand{\beginquestion}{\addtocounter{equation}{1}\begin{question}}
\newcommand{\beginsassumptions}{\addtocounter{equation}{1}\begin{sassumptions}}
\newcommand{\beginsassumption}{\addtocounter{equation}{1}\begin{sassumption}}
\newcommand{\beginconjecture}{\addtocounter{equation}{1}\begin{conjecture}}

\begin{document}

\title{A brief introduction to $p$-adic numbers}

\author{Stephen Semmes}

\date{}

\maketitle

\begin{abstract}
In this short survey we look at a few basic features of $p$-adic
numbers, somewhat with the point of view of a classical analyst.
In particular, with $p$-adic numbers one has arithmetic operations
and a norm, just as for real or complex numbers.
\end{abstract}

	Let ${\bf Z}$ denote the integers, ${\bf Q}$ denote the
rational numbers, ${\bf R}$ denote the real numbers, and ${\bf C}$
denote the complex numbers.  Also let $|\cdot |$ denote the usual
absolute value function or modulus on the complex numbers.

	On the rational numbers there are other absolute value functions
that one can consider.  Namely, if $p$ is a prime number, define the
\emph{$p$-adic absolute value function} $|\cdot |_p$ on ${\bf Q}$ by
$|x|_p = 0$ when $x = 0$, $|x|_p = p^{-k}$ when $x = p^k m/n$, where
$k$ is an integer and $m$, $n$ are nonzero integers which are not
divisible by $p$.  One can check that
\begin{equation}
	|x y|_p = |x|_p \, |y|_p
\end{equation}
and
\begin{equation}
	|x + y|_p \le |x|_p + |y|_p
\end{equation}
for all $x, y \in {\bf Q}$, just as for the usual absolute value function.
In fact, 
\begin{equation}
	|x + y|_p \le \max(|x|_p, |y|_p)
\end{equation}
for all $x, y \in {\bf Q}$.  This is called the \emph{ultrametric}
version of the triangle inequality.

	Just as the usual absolute value function leads to the
distance function $|x - y|$, the $p$-adic absolute value function
leads to the $p$-adic distance function $|x - y|_p$ on ${\bf Q}$.
With respect to this distance function, the rationals are not complete
as a metric space, and one can complete the rationals to get a larger
space ${\bf Q}_p$.  This is analogous to obtaining the real numbers by
completing the rationals with respect to the standard absolute value
function.  By standard reasoning the arithmetic operations and
$p$-adic absolute value function extend from ${\bf Q}$ to ${\bf Q}_p$,
with much the same properties as before.  In this manner one gets the
field of $p$-adic numbers.  As a metric space, ${\bf Q}_p$ is complete
by construction, and one can also show that closed and bounded subsets
of ${\bf Q}_p$ are compact.  This is also similar to the real numbers.

	Note that the set ${\bf Z}$ of integers forms a bounded subset
of ${\bf Q}_p$, in contrast to being an unbounded subset of ${\bf R}$.
In fact, each integer has $p$-adic absolute value less than or equal
to $1$.  There are general results about absolute value functions on
fields to the effect that if the absolute values of integers are bounded,
then they are less than or equal to $1$, and the absolute value function
satisfies the ultrametric version of the triangle inequality.  See p28-9
of \cite{Gouvea}.  In this case the absolute value function is said to
be \emph{non-Archimedian}.  If the absolute values of integers are not
bounded, as in the case of the usual absolute value function, then the
absolute value function is said to be \emph{Archimedian}.

	A related point is that the set ${\bf Z}$ of integers is a
\emph{discrete} subset of the real numbers.  It has no limit points,
and in fact the distance between two distinct integers is always at
least $1$.  This is not the case in ${\bf Q}_p$, where ${\bf Z}$ is
bounded, and hence precompact.  Now consider ${\bf Z}[1/p]$, the set
of rational numbers of the form $p^k n$, where $k$ and $n$ are
integers.  As a subset of ${\bf R}$, this is unbounded, and it also
contains nontrivial sequences which converge to $0$.  Similarly, as a
subset of ${\bf Q}_p$, it is unbounded and contains nontrivial
sequences which converge to $0$.  As a subset of ${\bf Q}_l$ when $l
\ne p$, ${\bf Z}[1/p]$ is bounded and hence precompact again.  Using
the diagonal mapping $x \mapsto (x, x)$, one can view ${\bf Z}[1/p]$
as a subset of the Cartesian product ${\bf R} \times {\bf Q}_p$.  In
this product, ${\bf Z}[1/p]$ is discrete again.  Indeed, if $a = p^k
b$ is a nonzero element of ${\bf Z}[1/p]$, where $k$, $b$ are integers
and $b$ is not divisible by $p$, then either $|a|_p \ge 1$, or $|a|_p
\le 1$, in which case $k \ge 0$, and $|a| \ge 1$.

	Similarly, $SL_n({\bf Z})$, the group of $n \times n$
invertible matrices with entries in ${\bf Z}$ and determinant $1$, is
a discrete subgroup of $SL_n({\bf R})$, the analogously-defined group
of matrices with real entries.  One can define $SL_n({\bf Z}[1/p])$
and $SL_n({\bf Q}_p)$ in the same manner, and using the diagonal
embedding $x \mapsto (x, x)$ again, $SL_n({\bf Z}[1/p])$ becomes a
discrete subgroup of the Cartesian product $SL_n({\bf R}) \times
SL_n({\bf Q}_p)$.

	There are fancier versions of these things for making ${\bf
Q}$ discrete, using ``adeles'', which involve $p$-adic numbers for all
primes $p$.  See \cite{Weil}.

	Now let us turn to some aspects of analysis.  With respect to
addition, ${\bf Q}_p$ is a locally compact abelian group, and thus has
a translation-invariant Haar measure, which is finite on compact sets,
strictly positive on nonempty open sets, and unique up to
multiplication by a positive real number.  As in \cite{Taibleson},
there is a rich Fourier analysis for real or complex-valued functions
on ${\bf Q}_p$, or ${\bf Q}_p^n$ when $n$ is a positive integer.

	Instead one can also be interested in ${\bf Q}_p$-valued
functions on ${\bf Q}_p$, or on a subset of ${\bf Q}_p$.  It is
especially interesting to consider functions defined by power series.
As is commonly mentioned, a basic difference between ${\bf Q}_p$ and
the real numbers is that an infinite series $\sum a_n$ converges if
and only if the sequence of terms $a_n$ tends to $0$ as $n$ tends to
infinity.  Indeed, the series converges if and only if the sequence of
partial sums forms a Cauchy sequence, and this implies that the terms
tend to $0$, just as in the case of real or complex numbers.  For
$p$-adic numbers, however, one can use the ultrametric version of the
triangle inequality to check that the partial sums form a Cauchy
sequence when the terms tend to $0$.  In particular, a power series
$\sum a_n \, x^n$ converges for some particular $x$ if and only if the
sequence of terms $a_n \, x_n$ tends to $0$, which is to say that
$|a_n|_p \, |x|_p^n$ tends to $0$ as a sequence of real numbers.

	Suppose that $\sum_{n = 0}^\infty a_n \, x^n$ is a power
series that converges for all $x$ in ${\bf Q}_p$, which is equivalent
to saying that $|a_n|_p r^n$ converges to $0$ as a sequence of real
numbers for all $r > 0$.  Thus we get a function $f(x)$ defined
on all of ${\bf Q}_p$, and we would like to make an analogy with
entire holomorphic functions of a single complex variable.  This
is somewhat like the situation of starting with a power series that
converges on all of ${\bf R}$, and deciding to interpret it as a 
function on the complex numbers instead.

	In fact, let us consider the simpler case of a power series
with only finitely many nonzero terms, which is to say a polynomial.
As in the case of complex numbers, it would be nice to be able to
factor polynomials.  The $p$-adic numbers ${\bf Q}_p$ are not
algebraically closed, and so in order to factor polynomials one 
can first pass to an algebraic closure.  It turns out that the $p$-adic
absolute value can be extended to the algebraic closure, while keeping
the basic properties of the absolute value.  See \cite{Cassels, Gouvea}.
The algebraic closure is not complete in the sense of metric spaces
with respect to the extended absolute value function, and one can
take a metric completion to get a larger field to which the absolute
values can be extended again.  A basic result is that this metric
completion is algebraically closed, so that one can stop here.
Let us write ${\bf C}_p$ for this new field, which is algebraically
closed and metrically complete.

	Once one goes to the algebraic closure, one can factor
polynomials.  On ${\bf C}_p$ one has this property and also one can
work with power series.  In particular, since the power series
$\sum_{n = 0}^\infty a_n \, x^n$ converges on all of ${\bf Q}_p$, it
also converges on all of ${\bf C}_p$, so that $f(x)$ can be extended
in a natural way to ${\bf C}_p$.  For that matter, one can start with
a power series that converges on all of ${\bf C}_p$, where the
coefficients are allowed to be in ${\bf C}_p$, and not just ${\bf Q}_p$.

	Under these conditions, the function $f$ can be written as
a product of an element of ${\bf C}_p$, factors which are equal to $x$,
and factors of the form $(1 - \lambda_j \, x)$, where the
$\lambda_j$'s are nonzero elements of ${\bf C}_p$.  In other words,
the factors of $x$ correspond to a zero of some order at the origin,
while the factors of the form $(1 - \lambda_j \, x)$ correspond
to zeros at the reciprocals of the $\lambda_j$'s.  If $f(x)$ is a
polynomial, then there are only finitely many factors, and this
statement is the same as saying that ${\bf C}_p$ is algebraically
closed.  In general, each zero of $f(x)$ is of finite order, and
there are only finitely many zeros within any ball of finite radius
in ${\bf C}_p$.  Thus the set of zeros is at most countable, and 
this condition permits one to show that the product of the factors
mentioned above converges when there are infinitely many factors.

	This representation theorem can be found on p113 of
\cite{Cassels} and on p209 of \cite{Gouvea}.  It is analogous to
classical results about entire holomorphic functions of a complex
variable, with some simplifications.  In the complex case, it is
necessary to make assumptions about the growth of an entire function
for many results, and the basic factors often need to be more
complicated in order to have convergence of the product.  See
\cite{Ahlfors, Veech} concerning entire holomorphic functions of
a complex variable.


\begin{thebibliography}{99}


\addcontentsline{toc}{section}{References}

\bibitem {Ahlfors} L.~Ahlfors, {\it Complex Analysis}, third edition,
McGraw-Hill, 1979.

\bibitem {Borel} A.~Borel, {\it Essays in the History of Lie Groups
and Algebraic Groups}, History of Mathematics {\bf 21}, American
Mathematical Society and London Mathematical Society, 2001.

\bibitem {Brown1} K.~Brown, {\it Buildings}, Springer Monographs in
Mathematics, Springer-Verlag, 1989.

\bibitem {Brown2} K.~Brown, {\it What is a building?}, Notices of the
American Mathematical Society {\bf 49} (2002), 1244-1245.

\bibitem {Cassels} J.~Cassels, {\it Local Fields}, London Mathematical
Society Student Texts {\bf 3}, Cambridge University Press, 1986.

\bibitem {Goldberg} R.~Goldberg, {\it Methods of Real Analysis},
Blaisdell, 1964.

\bibitem {Gouvea} F.~Gouv\^ea, {\it $p$-Adic Numbers: An Introduction},
Universitext, Springer-Verlag, 1993.

\bibitem {Rudin} W.~Rudin, {\it Principles of Mathematical Analysis},
third edition, McGraw-Hill, 1976.

\bibitem {SS} S.~Semmes, {\it Some topics pertaining to algebras of
linear operators}, arXiv.org archive, math.CA/0211171.

\bibitem {Serre1} J.-P.~Serre, {\it Corps Locaux}, Hermann, 1962;
{\it Local Fields}, translated from the French by M.~Greenberg, Graduate
Texts in Mathematics {\bf 67}, Springer-Verlag, 1979.

\bibitem {Serre2} J.-P.~Serre, {\it Cours d'Arithm\'etique}, Presses
Universitaires de France, 1970; {\it A Course in Arithmetic}, Graduate
Texts in Mathematics {\bf 7}, 1973.

\bibitem {Serre3} J.-P.~Serre, {\it Cohomologie des groupes discrets},
in {\it Prospects in Mathematics}, 77--169, Annals of Mathematics
Studies {\bf 70}, Princeton University Press, 1971.

\bibitem {Taibleson} M.~Taibleson, {\it Fourier Analysis on Local Fields},
Mathematical Notes {\bf 15}, Princeton University Press, 1975.

\bibitem {Veech} W.~Veech, {\it A Second Course in Complex Analysis},
Benjamin, 1967.

\bibitem {Weil} A.~Weil, {\it Basic Number Theory}, second edition,
Springer-Verlag, 1973.


\end{thebibliography}
\end{document}